\documentclass{ifacconf}
\usepackage{amsfonts}
\usepackage{enumitem}
\usepackage{amsmath}
\usepackage{amssymb}
\usepackage{algorithm}
\usepackage{booktabs}
\usepackage{graphicx,color}
\usepackage{natbib}

\newtheorem{theorem}{Theorem}

\newtheorem{remark}{Remark}
\newtheorem{assumption}{Assumption}

\newcommand{\ms}{{{\rm\bf q}}}

\newcommand{\solvelp}{{\texttt{Solve\textsubscript{LP}}}}

\renewcommand{\Pr}{{\mathbb{P}}}
\newcommand{\Bad}{{\textrm{Bad}}}

\renewcommand{\natural}{{\mathbb{N}}}

\newcommand{\until}[1]{\{1,\ldots,#1\}}

\newcommand{\EE}{\mathcal{E}}

\newcommand{\VV}{\mathcal{V}}

\renewcommand{\theenumi}{(\roman{enumi})} 

\begin{document}
\begin{frontmatter}

% \title{\MC{Tentative: Randomized Constraints Consensus for Distributed Linear Programming in the Presence of Uncertainty}
% \thanksref{footnoteinfo}} 
% \title{Randomized Constraints Consensus for Distributed Uncertain Linear Programs\thanksref{footnoteinfo}} 
\title{Randomized Constraints Consensus for Distributed Robust Linear Programming\thanksref{footnoteinfo}} 


\thanks[footnoteinfo]{This work is supported by the European Research Council
  (ERC) under the European Union's Horizon 2020 research and innovation
  programme, grant agreement No 638992 - OPT4SMART, (GN)  and by a grant from
  the Singapore National Research Foundation (NRF) under the ASPIRE project,
  grant No NCR-NCR001-040 (MC\&RB). % \\
  % \emph{Email addresses:} \texttt{Chamanbaz@sutd.edu.sg} (Mohammadreza
  % Chamanbaz) and \texttt{giuseppe.notarstefano@unisalento.it} (Giuseppe
  % Notarstefano) and \texttt{bouffanais@sutd.edu.sg} (Roland Bouffanais).
}

\author[ARAKUT,SUTD]{Mohammadreza Chamanbaz} 
\author[unisalento]{Giuseppe Notarstefano} 
\author[SUTD]{Roland Bouffanais}

\address[ARAKUT]{Arak University of Technology,
  Arak, Iran\\ \texttt{Chamanbaz@arakut.ac.ir}}
\address[unisalento]{Department of Engineering, Universit\`{a} del Salento, Lecce, Italy\\\texttt{giuseppe.notarstefano@unisalento.it}}
\address[SUTD]{Singapore University of Technology and Design,
  Singapore\\ \texttt{\{Chamanbaz,bouffanais\}@sutd.edu.sg} }

%-----------------------------------------------------------------------------------------

\begin{abstract}                % Abstract of not more than 250 words.
  In this paper we consider a network of processors aiming at cooperatively
  solving linear programming problems subject to uncertainty. Each node only
  knows a common cost function and its local uncertain constraint set. We
  propose a randomized, distributed algorithm working under time-varying,
  asynchronous and directed communication topology. The algorithm is based on a
  local computation and communication paradigm. At each communication round,
  nodes perform two updates: (i) a verification in which they check---in a
  randomized setup---the robust feasibility (and hence optimality) of the
  candidate optimal point, and (ii) an optimization step in which they exchange
  their candidate bases (minimal sets of active constraints) with neighbors and
  locally solve an optimization problem whose constraint set includes: a sampled
  constraint violating the candidate optimal point (if it exists), agent's
  current basis and the collection of neighbor's basis. As main result, we show
  that if a processor successfully performs the verification step for a
  sufficient number of communication rounds, it can stop the algorithm since a
  consensus has been reached. The common solution is---with high
  confidence---feasible (and hence optimal) for the entire set of uncertainty
  except a subset having arbitrary small probability measure.  We show the
  effectiveness of the proposed distributed algorithm on a multi-core platform
  in which the nodes communicate asynchronously.
\end{abstract}


\end{frontmatter}

%-----------------------------------------------------------------------------------------

\section{Introduction}

Robust optimization plays an important role in several areas such as
estimation and control and has been widely investigated. Its rich literature
dates back to the $1950$s, see~\cite{ben2009robust} and references therein. Very
recently, there has been a renewed interest in this topic in a parallel and/or
distributed framework.
% However, there are only few papers considering a distributed framework for
% solving robust optimization problems.
In \cite{lee2013distributed}, a synchronous distributed
random projection algorithm with almost sure convergence is proposed for the
case where each node has independent cost function and (uncertain)
constraint. Since the distributed algorithm relies on extracting random samples
from an uncertain constraint set, several assumptions on random set, network
structure and agent weights are made to prove almost sure convergence. The
synchronization of update rule relies on a central clock to coordinate the step
size selection. To circumvent this limitation the same authors in
\cite{lee2016asynchronous} present an asynchronous random projection algorithm
in which a gossip-based protocol is used to desynchronize the step size
selection. The proposed algorithms in
\citep{lee2013distributed,lee2016asynchronous}, require computing projection
onto the constraint set at each iteration which is computationally demanding if
the constraint set does not have a simple structure such as half space or
polyhedron.
In \cite{carlone2014distributed}, a
parallel/distributed scheme is considered for solving an uncertain optimization
problem by means of the scenario approach
\citep{calafiore_uncertain_2004}. The scheme consists of extracting a number of
samples from the uncertain set and assigning them to nodes in a network. Each
node is assigned a portion of the extracted samples. Then, a variant of the
constraints consensus algorithm introduced in \cite{notarstefano2011distributed} is used
to solve the deterministic optimization problem. 
%
% In this approach, each node of the network is
% given a subset of sampled optimization problem. Nodes exchange their set of
% active constraints and solve their local optimization problem till all nodes
% reach a consensus on the set of active constraints. It is clear that the number
% of constraints in the optimization problem each node has to solve increases with
% the iteration counter giving rise to a more complex optimization problem. In
% fact, at the extreme case, all the samples can be active and hence, each node
% has to solve an optimization problem with the same complexity as the original
% scenario optimization problem.
%
A similar parallel framework for solving convex optimization problems with one
uncertain constraint via the scenario approach is considered in
\cite{you_networked_2016}. In this setup, the sampled optimization problem is
solved in a distributed way by using a primal-dual subgradient (resp. random
projection) algorithm over an undirected (resp. directed) graph. We remark that in
\cite{carlone2014distributed,you_networked_2016}, constraints and cost function
of all agents are identical. In \cite{burger2014polyhedral}, a cutting plane
consensus algorithm is introduced for solving convex optimization problem where
constraints are distributed to the network processors and all processors have
common cost function. In the case where constraints are uncertain, a worst-case
approach based on a pessimizing oracle is used. The oracle relies on the
assumption that constraints are concave with respect to uncertainty vector and
the uncertainty set is convex.
A distributed scheme based on the scenario approach is introduced in
\cite{margellos2016distributed} in which random samples are extracted by each
node from its local uncertain constraint set and a distributed proximal
minimization algorithm is designed to solve the sampled optimization
problem. 
The number of samples required to guarantee robustness can be large if the
probabilistic levels defining robustness of the solution---accuracy and
confidence levels---are stringent, possibly leading to a computationally demanding
sampled optimization problem at each node.

% CONTRIBUTIONS 

% -) we sample locally on local constraint
% -) we don't need to sample everything at the beginning
% -) we don't need assumptions as Nedic on the uncertain set
% -) we use most of the samples only for verification (computationally inexpensive, while few samples in the
%    local optimization

The main contribution of this paper is the design of a fully distributed
algorithm to solve an uncertain linear program in a network with directed and
asynchronous communication. The problem under investigation is a linear program
in which the constraint set is the intersection of local uncertain constraints,
each one known only by a single node. Starting from a deterministic constraint
exchange idea introduced in \cite{notarstefano2011distributed}, the algorithm
proposed in this paper introduces a randomized, sequential approach in which
each node: (i) locally performs a probabilistic verification step (based on a
local sampling of its uncertain constraint set), and (ii) solves a local,
deterministic optimization problem with a limited number of constraints. If
suitable termination conditions are satisfied, we are able to prove that the
nodes agree on a common solution which is probabilistically feasible and optimal
with high confidence.
As compared to the literature above, the proposed algorithm has three main
advantages. First, no assumptions are needed on the probabilistic nature of the
local constraint sets. Second, each node can sample locally its own uncertain
set. Thus, no central unit is needed to extract samples and no common constraint
set needs to be known by the nodes. Third and final, nodes do not need to
perform the whole sampling at the beginning and subsequently solve the (deterministic)
optimization problem. Online extracted samples are used only for
verification, which is computationally inexpensive. The optimization is
performed always on a number of constraints that remains constant at each node
and depends only on the dimension of the decision variable and on the number of
node neighbors.

% our algorithm can use any LP solver locally (inherited from the constraints consensus
% \textbf{Paper Organization}\\
The paper is organized as follows. In Section~\ref{sec: problem formulation}, we
formulate the uncertain distributed linear program (LP). Section~\ref{sec:
  distributed algorithm} presents our distributed sequential randomized
algorithm for finding a solution---with probabilistic robustness---to the
uncertain distributed LP. The probabilistic convergence properties of the
distributed algorithm are investigated in Section~\ref{sec: analysis of the
  algorithm}. Finally, extensive numerical simulations are performed in
Section~\ref{sec: numerical simulations} to prove the effectiveness of the
proposed methodology.
% Some concluding remarks are collected in Section~\ref{sec: conclusions}.

%-----------------------------------------------------------------------------------------

\section{Problem Formulation}\label{sec: problem formulation}
We consider a network of processors 
% characterized by directed communication graph
% $\mathcal{G}=\{\mathcal{V},\mathcal{E}\}$ where
% $\mathcal{V}\doteq\{1,\ldots,n\}$ is the set of nodes (processors) and
% $\mathcal{E}\subseteq\mathcal{V}\times\mathcal{V}$ is the set of edges
% (communication links). Processors have
with limited computation and/or communication capabilities that aim at
cooperatively solving the following uncertain linear program
\begin{align}\nonumber
\underset{\theta}{\min}&~~c^T\theta\\ \label{eq: uncertain LP}
\text{subject to}&~~A_i^T(q)\theta\leq b_i(q), \; \forall q\in\mathbb{Q},~i\in\{1,\ldots,n\},
\end{align}
where $\theta\in\Theta\subset\mathbb{R}^{d}$ is the vector of decision
variables, $q\in\mathbb{Q}$ is the uncertainty vector, $c\in\mathbb{R}^d$
defines the objective direction, $A_i(q)\in\mathbb{R}^{m_i\times d}$ and
$b_i(q)\in\mathbb{R}^{m_i}$, with $m_i\geq d$, define the (uncertain) constraint
set of agent $i\in\until{n}$. 
% Constraint set of agent $i$ contains $m_i\geq d$ inequality constraints.
Processor $i$ has only knowledge of a constraint set defined by $A_i(q)$ and
$b_i(q)$ and the objective direction $c$ (which is the same for all nodes).
Each node runs a local algorithm and by exchanging limited information with
neighbors, all nodes converge to the same solution.
We want to stress that there is no (central) node having access to all
constraints.
% or that nodes cannot extract a common uncertainty vector
%$q\in\mathbb{Q}$.
%
We make the following assumption regarding
problem~\eqref{eq: uncertain LP}.
\begin{assumption}[Non-degeneracy]~\label{assum: nondegeneracy} The minimum
  point of any subproblem of~\eqref{eq: uncertain LP} with at least $d$
  constraints is unique and there exist only $d$ constraints intersecting at the
  minimum point.
\end{assumption}
%{
% Assumption~\ref{assum: nondegeneracy} can be relaxed by introducing lexicographical order. \textbf{Maybe we can explain a bit what lexicographical
% order is??}}
We let the nodes communicate according to a time-dependent, directed
communication graph 
% If communication links do not change over time, the communication graph $\mathcal{G}$ is called \emph{static}, otherwise it is called \emph{time-varying}. The time-varying communication case is modeled by a time-varying directed graph
$\mathcal{G}(t)=\{\mathcal{V},\mathcal{E}(t)\}$ where $t\in\natural$ is a
universal time, $\VV=\until{n}$ is the set of agent identifiers and
$(i,j)\in\EE(t)$ indicates that $i$ send information to $j$ at time $t$. The
time-varying set of incoming (resp. outgoing) neighbors of node $i$ at time $t$,
$\mathcal{N}_\text{in}(i,t)$ ($\mathcal{N}_\text{out}(i,t)$), is defined as the
set of nodes from (resp. to) which agent $i$ receives (resp. transmits) information at time
$t$.
% The set of incoming (outgoing) neighbors is denoted by
% $\mathcal{N}_\text{in}(i,t)$ ($\mathcal{N}_\text{out}(i,t)$).
%In a static graph, the \emph{distance} between nodes $i$ and $j$---the minimum number of edges between node $i$ and $j$---is denoted by $\text{dist}(i,j)$. The maximum distance over all nodes is called \emph{diameter} of the graph and is denoted by $\text{diam}(\mathcal{G})$. 
A directed static graph is said to be \emph{strongly connected} if there exists
a directed path (of consecutive edges) between any pair of nodes in the
graph. For time-varying graphs we use the notion of \emph{uniform joint strong
  connectivity} formally defined next.

\begin{assumption}[Uniform joint strong connectivity]~\\
\label{assum:graph}
There exists an integer $L\ge1$ such that the graph
%  $\bigg(\mathcal{V}, \bigcup_{\tau=tL}^{(t+1)L-1} \mathcal{E}(\tau)\bigg)$ 
$\bigg(\mathcal{V}, \bigcup_{\tau=t}^{t+L-1} \mathcal{E}(\tau)\bigg)$ is strongly
  connected for all $t\ge0$.
\end{assumption}

There is no assumption on how uncertainty $q$ enters problem \eqref{eq:
  uncertain LP} making it computationally difficult to solve. In fact, if the
uncertainty set $\mathbb{Q}$ is an uncountable set, problem \eqref{eq: uncertain
  LP} is a semi-infinite optimization problem involving infinite number of
constraints. In general, there are two main paradigms to solve an uncertain
optimization problem of form \eqref{eq: uncertain LP}. The first approach is a
deterministic worst-case paradigm in which the constraints are enforced to hold
for \emph{all} possible uncertain parameters in the set $\mathbb{Q}$. This
approach is computationally intractable for cases where uncertainty does not
appear in a ``simple'' form, e.g. affine, multi-affine, convex, etc.  The second
approach is a probabilistic approach where uncertain parameters are considered
to be random variables and the constraints are enforced to hold for the entire
set of uncertainty except a subset having an arbitrary small probability
measure. In this paper, we follow a probabilistic approach and present a
distributed tractable randomized setup for finding a solution---with desired
probabilistic properties---for the optimization problem \eqref{eq: uncertain
  LP}.

\noindent
{\bf Notation \\}
% \begin{remark}[Notations]
The constraint set of  agent $i$ is defined by %$H^i(q)\in\mathbb{R}^{m_i\times(d+1)}$
\[
H^i(q)\doteq[A_i(q),b_i(q)].
\]
Throughout this paper, we use capital italic letter,
e.g. $H^i(q)\doteq[A_i(q),b_i(q)]$ to denote a collection of half spaces and
capital calligraphic letter,  $\mathcal{H}^i(q)$ to denote the set induced by
half spaces, i.e. $\mathcal{H}^i(q)\doteq\{\theta\in\mathbb{R}^d:A_i(q)\leq
b_i(q)\}$. We note that, with this notation, if $A=B\cup C$ with $B$ and $C$
being collection of half spaces,  then $\mathcal{A}=\mathcal{B}\cap
\mathcal{C}$, that is, the set induced by the union of constraint sets $B$ and
$C$ is the intersection of $\mathcal{B}$ and $\mathcal{C}$.
Finally $J(H)$ is the smallest value of $c^T\theta$ while
$\theta\in\mathcal{H}$. 
The linear program specific to each agent $i\in\mathcal{V}$ is fully
characterized by the pair $(H^i(q),c)$ (note that $c$ defines the objective
direction which is the same for all nodes).

\section{Randomized Constraints Consensus}\label{sec: distributed algorithm}
In this section, we present a distributed, randomized algorithm for solving the uncertain
linear program (LP) \eqref{eq: uncertain LP} in a probabilistic sense. 
% The algorithm is executed locally on each processor till all nodes converge to a
% common solution with desired probabilistic properties.
First, recall that the solution of a linear program of the form \eqref{eq:
  uncertain LP} can be identified by at most $d$ active constraints ($d$ being
the dimension of the decision variable).
This concept is formally characterized by the notion of \emph{basis}.  Given a
collection of constraints $H$, a subset $B\subseteq H$ is a basis of $H$ if the
optimal cost of the LP problem defined by $(H,c)$ is identical to the one
defined by $(B,c)$, and the optimal cost decreases if any constraint is removed
from $B$.
% We remark that under Assumption~\ref{assum: nondegeneracy}, the set of
% all possible active constraints coincides with the basis.
We define a primitive $[\theta^*,B]=\solvelp(H,c)$ which solves the LP problem
defined by the pair $(H,c)$ and returns back the optimal point $\theta^*$ and
the corresponding basis $B$.

Note that, since the uncertainty set is uncountable, it is in general very
difficult to verify if a candidate solution is feasible for the entire set of
uncertainty or not.
We instead use a randomized approach based, on Monte Carlo simulation, to check
probabilistic feasibility. %  of the candidate solution.
The distributed algorithm we propose has a probabilistic nature consisting of
two main steps: verification and optimization. The main idea is the following. A
node has a candidate basis and candidate solution point. First, it verifies if
the candidate solution point belongs to its local uncertain set with high
probability. Then, it collects bases from neighbors and solves an LP with its
basis and its neighbors' bases as constraint set. If the verification step was
not successful, the first violating constraint is also added to the problem.

% In the verification step, robust feasibility of a candidate solution
% $\theta^i(t)$ is checked. That is, node $i$ checks---with high
% probability---if $\theta^i(t)\in\mathcal{H}^i(q)$ for all $q\in\mathbb{Q}$.
%

Formally, we assume that $q$ is a random variable and a probability measure
$\Pr$ over the Borel $\sigma-$algebra of $\mathbb{Q}$ is given.
In the verification step each agent $i$ generates $M_{k_i}$ independent and
identically distributed (i.i.d) random samples from the set of uncertainty
\[
\ms_{k_i}\doteq\{q^{(1)},\ldots,q^{(M_{k_i})}\}\in\mathbb{Q}^{M_{k_i}},
\] 
according to the measure $\Pr$, where ${k_i}$ is a local counter keeping track
of the number of times the verification step is performed and
$\mathbb{Q}^{M_{k_i}}\doteq\mathbb{Q}\times\ldots\times\mathbb{Q}$ ($M_{k_i}$
times).  Using a Monte Carlo algorithm, node $i$ checks feasibility of the
candidate solution $\theta^i(t)$ only at the extracted samples. If a violation
happens, the first violating sample is used as a \emph{violation certificate}.
In the optimization step, agent $i$ transmits its current basis to all outgoing
neighbors and receives bases from incoming ones. Then, it solves an LP problem
whose constraint set is composed of: \emph{i)} a constraint constructed at the
violation certificate (if it exists) \emph{ii)} its current basis and
\emph{iii)} the collection of bases from all incoming neighbors.
Node $i$ repeats these two steps until a termination condition is satisfied,
namely if the candidate basis has not changed for $2 n L+1$ times, with $L$
defined in Assumption~\ref{assum:graph}.
% These two steps, i.e. verification and optimization
% are repeated till all nodes converge to a common solution which is robustly
% feasible---with desired (high) probability---for all constraints across
% different nodes. 
%
The distributed algorithm is formally presented in Algorithm~\ref{alg:
  distributed algorithm}.
\begin{algorithm}[h]
\caption{Randomized Constraints Consensus}
\label{alg: distributed algorithm}
\textbf{Input:}{ $(H^i(q),c),\varepsilon_i,\delta_i$}

\textbf{Output:}{ $\theta_\text{sol}$}

\textbf{Initialization:\\}
Set $ k_i=1,  [\theta^i(1),B^i(1)] = \solvelp(H^i(\mathbf{0}),c)$

\textbf{Evolution:}
\begin{enumerate}[leftmargin=*]
\item
\textbf{Verification:}\label{step :verification}

\begin{itemize}[leftmargin=*]
\item
If $\theta^i(t)=\theta^i(t-1)$, set $q_{\texttt{viol}}=\emptyset$ and goto~\ref{step: optimization}
\item 
Extract 
\begin{equation}\label{eq:sample bound Mk}
M_{k_i}\geq \frac{2.3+1.1\ln k_i+\ln \frac{1}{\delta_i}}{\ln \frac{1}{1-\varepsilon_i}}
\end{equation}
i.i.d samples $\ms_{k_i} = \{q_{k_i}^{(1)},\ldots,q_{k_i}^{(M_{k_i})}\}$
\item
If $\theta^i(t)\in\mathcal{H}^i(q_{k_i}^{(\ell)})$ for all
$\ell=1,\ldots,M_{k_i}$, set $q_{\texttt{viol}}=\emptyset$; else, set
$q_{\texttt{viol}}$ as the first sample for which $\theta^i(t)\notin
\mathcal{H}^i(q_{\texttt{viol}})$
\item
Set $k_i=k_i+1$
\end{itemize}
\item
\textbf{Optimization:} \label{step: optimization}
\begin{itemize}[leftmargin=*]
\item 
Transmit $B^i(t)$ to $j\in\mathcal{N}_\text{out}(i,t)$ and acquire incoming neighbors basis $Y^i(t)\doteq\cup_{j\in\mathcal{N}_\text{in}(i,t)}B^j$
\item
$[\theta^i(t+1),B^i(t+1)]=$\\
$~~~~~~~~~~~~~\solvelp(H^i(q_{\texttt{viol}})\cup B^i(t)\cup Y^i(t),c)$
%\item 
%Set $B^i(t+1)$ as the set of minimal basis of $h^i(q_{\texttt{viol}})\cup B^i(t)\cup Y^i(t)$. 
\item
If  $\theta^i(t+1)$ has not changed for $2nL+1$ times and $q_\texttt{viol}=\emptyset$, return  $\theta_\text{sol}=\theta^i(t+1)$
\end{itemize}
%\item \label{step: update}
%\textbf{Update:} Set $t=t+1$ and goto~\ref{step :verification}.
\end{enumerate}
\end{algorithm} 
The counter $k_i$ counts the number of times the verification step is called.
We remark that if at some $t$ the candidate solution has not changed, that is
$\theta^i(t)=\theta^i(t-1)$, then $\theta^i(t-1)$ has successfully satisfied the
verification step and $q_\texttt{viol}=\emptyset$ at time $t-1$ and therefore
there is no need to check it again.  

% \MC{It is clear that, if $q_\texttt{viol}\neq\emptyset$, then the constraint formed
% at $q_\texttt{viol}$, $H^i(q_\texttt{viol})$ , which is infeasible for
% $\theta^i(t-1)$, will change the solution of the optimization problem in
% step~\ref{step: optimization} of the algorithm and hence
% $\theta^i(t)\neq\theta^i(t+1)$.
% %and since it is not changed, there is no need to verify it again in a Monte Carlo simulation.  
% Therefore, in the  step~\ref{step :verification} of the distributed algorithm
% its robustness is checked in the Monte Carlo simulation only if the candidate
% solution has changed at the previous step. 
% }

%\begin{remark}[Computational complexity]
%  % It is computationally inexpensive to verify the feasibility of a candidate
%  % solution $\theta^i(t)$ for a large number of samples extracted from the
%  % local uncertain constraint set.
%  The verification step of
%  Algorithm~\ref{alg: distributed algorithm} is not computationally
%  demanding even though the number of extracted samples is large, since it
%  only involves evaluating extracted inequalities. 
%  %
%  The optimization step of Algorithm~\ref{alg: distributed
%    algorithm}, even though more demanding, contains a fixed number of
%  constraints, which only depends on the number of node neighbors and dimension
%  of decision variables $d$. In particular, it does not depend on $k_i$.
%\end{remark}

\begin{remark}[Asynchronicity]
  The distributed algorithm presented in this section is completely
  asynchronous. Indeed, time $t$ is just a universal time that does not need to
  be known by the nodes. The time-dependent jointly connected graph then captures
  the fact that nodes can perform computation at different speeds.
\end{remark}

\begin{remark}
  In the deterministic constraints consensus algorithm presented in
  \cite{notarstefano2011distributed}, at each iteration of the algorithm, the
  original constraint set of the node needs to be taken into account in the
  local optimization problem. Here, we can drop this requirement because of the
  verification step.
\end{remark}

%-----------------------------------------------------------------------------------------

\section{Analysis of Randomized Constraints Consensus algorithm} \label{sec: analysis of the algorithm}
In this section, we analyze the convergence properties of the distributed
algorithm and investigate the probabilistic properties of the solution computed
by the algorithm.  

%\begin{definition}[Constraints violation probability]
%The constraints violation probability of the basis $B$ for the constraint set $H(q)$ is defined as 
%\begin{equation}\label{eq:constraint violation}
%V_c(B,H(q))\doteq\Pr\{q\in\mathbb{Q}:\mathcal{B}\not\subset\mathcal{H}(q)\}.
%\end{equation}
%\end{definition}
%
%\begin{definition}[Cost violation probability]
%The cost violation probability of the basis $B$ for the constraint set $H(q)$ is defined as 
%\begin{equation}\label{eq:cost violation}
%V(B,H(q))\doteq\Pr\{q\in\mathbb{Q}:J(B\cup H(q))>J(B)\}.
%\end{equation}
%\end{definition}

\begin{theorem}\label{thm: convergence }
Let Assumptions~\ref{assum: nondegeneracy} and~\ref{assum:graph} hold. Given the
probabilistic levels $\varepsilon_i>0$ and $\delta_i>0$, $i = 1,\ldots,n$, let
$\varepsilon= \sum_{i=1}^n\varepsilon_i$ and $\delta= \sum_{i=1}^n\delta_i$.
Then, the following statements hold
\begin{enumerate}
\item
Along the evolution of Algorithm~\ref{alg: distributed algorithm}, the cost
$J(B^i(t))$ at each node $i\in\{1,\ldots,n\}$ is monotonically
non-decreasing. That is, $J(B^i(t+1))\geq J(B^i(t))$. 
\item 
The  cost $J(B^i(t))$ for all $i\in\{1,\ldots,n\}$  converges to a common value
asymptotically. That is,  $\lim_{t\rightarrow\infty}J(B^i(t))=\bar{J}$ for all
$i\in\{1,\ldots,n\}$. 
%\item
%%\MC{If we include non-degeneracy assumption we can have this statement. \\}
%The  basis and candidate solution of all nodes converge to a common value asymptotically. That is $B^1(t)=\ldots=B^n(t)$ and $\theta^1(t)=\ldots=\theta^n(t)$ as $t\rightarrow\infty$.
\item
If the candidate solution  of node $i$, $\theta^i(t)$, has not changed for
$2Ln+1$ communication rounds, all nodes have a common candidate
solution $\theta_\text{sol}$\footnote{We remark that this value becomes
  $2\times$(graph diameter)+1
  for fixed graphs.}. 
\item 
The following inequality holds for $\theta_\text{sol}$
\begin{align*}
\Pr^M\bigg\{&\ms\in\mathbb{Q}^M:\\
&\Pr\bigg\{q\in\mathbb{Q}:\theta_\text{sol}\notin \bigcap_{i=1}^n\mathcal{H}^i(q)\bigg\}\leq\varepsilon \bigg\}
\geq1-\delta,
\end{align*}
where $M$ is the cardinality of the collection of multisamples of all agents. 

\item
Let $B_\text{sol}$ be the basis corresponding to $\theta_\text{sol}$. 
The following inequality holds for $B_\text{sol}$
%\MC{\\ I am not sure how to state this using the point $\theta_\text{sol}$}
\begin{align*}
\Pr^M\bigg\{\ms\in\mathbb{Q}^M:\Pr\bigg\{q\in\mathbb{Q}:J&\big(B_\text{sol}\cup H(q)\big)>J(B_\text{sol})\bigg\}\\
&\leq\varepsilon \bigg\}\geq1-\delta
\end{align*}
where $H(q)\doteq \bigcup_{i=1}^n H^i(q)$ and $M$ is the cardinality of the collection of multisamples of all agents.

\end{enumerate}
\end{theorem}

%\textbf{Proof.\\}~ 
\textbf{Proof of first statement:\\} The set of constraints at time $t+1$
consists of the node current basis $B^i(t)$, the collection of neighbors' bases
$Y^i(t)\doteq\cup_{j\in\mathcal{N}_\text{in}(i,t)}B^j$ and
$H^i(q_\texttt{viol})$. Since the basis at time $t$, $B^i(t)$, is part of
the constraint set for computing $B^i(t+1)$, $J(B^i(t+1))$ cannot be smaller than
$J(B^i(t))$.

\textbf{Proof of second and third statements:\\}
%\MC{Please check this\\}
%if $\ell_1\in\mathcal{N}_\text{out}(i,t_0)$ then  $J(B^i(t_0))\leq J(B^{\ell_1}(t_0+1))$ because the constraint set of node $\ell_1$ at time $t_0+1$ is a superset of constraints set of node $i$ at time $t$.
Since the graph is uniformly jointly strongly connected, for any pairs of nodes
$u$ and $v$ and for any $t>0$, there exists a time-dependent path from $u$ to
$v$ \citep{hendrickx2008graphs}---a sequence of nodes $\ell_1,\ldots,\ell_k$  and a sequence of time
instances $t_1,\ldots,t_{k+1}$ with $t\leq t_1<\ldots<t_{k+1}$, such
that the directed edges
$\{(u,\ell_1),(\ell_1,\ell_2),\ldots,(\ell_k,v)\}$ belongs to the
directed graph at time instances $\{t_1,\ldots,t_{k+1}\}$, respectively---of length at most $nL$.
We recall that $n$ is the number of nodes and $L$ is defined in
Assumption~\ref{assum:graph}. Consider nodes $i$ and $p$. If
$\ell_{1}\in\mathcal{N}_\text{out}(i,t_0)$, then
$J(B^i(t_0))\leq J(B^{\ell_1}(t_0+1))$ as the constraint set of node
$\ell_{1}$ at time $t_0+1$ is a superset of the constraint set of node $i$
at time $t_0$. Iterating this argument, we obtain
$J(B^i(t_0))\leq J(B^p(t_0+nL))$.
%Now let $\ell_2\in\mathcal{N}_\text{out}(p,t_0+nL)$ hence, $J(B^i(t_0))\leq J(B^p(t_0+nL)\leq J(B^{\ell_2}(t_0+nL+1))  $. 
Again since the graph is uniformly jointly strongly connected, there will be a time varying path of length at most $nL$ from node $p$ to node $i$. Therefore,
%iterating this $\text{dist}(p,i)$ times, we arrive at $J(B^i(t_0))\leq J(B^p(t_0+\text{dist}(i,p)))\leq J(B^i(t_0+\text{dist}(i,p)+\text{dist}(p,i))) $. .
%Iterating this argument $\text{dist}(i,p)$ times \MC{is distance defined for a time varying graphs??}, we get $J(B^i(t_0))\leq J(B^p(t_0+\text{dist}(i,p)))$. Now let $\ell_2\in\mathcal{N}_\text{out}(p,t_0+\text{dist}(i,p))$ hence, $J(B^i(t_0))\leq J(B^p(t_0+\text{dist}(i,p)))\leq J(B^{\ell_2}(t_0+\text{dist}(i,p)+1))  $. Again iterating this $\text{dist}(p,i)$ times, we arrive at $J(B^i(t_0))\leq J(B^p(t_0+\text{dist}(i,p)))\leq J(B^i(t_0+\text{dist}(i,p)+\text{dist}(p,i))) $. Since the graph is uniformly jointly strongly connected, for any pairs of nodes $u$ and $v$ and for any $t>0$, there exists a time-dependent path from $u$ to $v$---a path $u,\ell_1,\ldots,\ell_k,v$ such that $(\ell_t,\ell_{t+1})\in\mathcal{E}(t+1)$---of length at most $nL$ \MC{cite the relevant paper}. We recall that $n$ is the number of nodes and $L$ is defines in Assumption \ref{assum:graph}. Therefore,
\[
J(B^i(t_0))\leq J(B^p(t_0+nL))\leq J(B^i(t_0+2nL)).
\]
Two scenarios can happen proving respectively statements (ii) and (iii).

%\begin{enumerate}
% \item[1.] 
If $J(B^i(t_0))\neq J(B^i(t_0+2nL))$, then $J(B^i(t_0)) < J(B^i(t_0+2nL))$
  which means the cost at node $i$ is strictly increasing.  Denote by $J^*$ the
  optimal cost associated to problem \eqref{eq: uncertain LP}. Since the problem
  at each node can be considered as a sub-problem of~\eqref{eq: uncertain LP},
  the sequence $\{J(B^i(t))\}_{t>0}$ is convergent and has a limit point
  $\bar{J}^i\leq J^*$, i.e.
  $\lim_{t\rightarrow\infty}J(B^i(t))\rightarrow \bar{J}^i$ for all
  $i\in\mathcal{V}$.  In what follows, we prove that the limit point of all
  agents are the same, that is $\bar{J}^1=\ldots=\bar{J}^n$. 
We follow a similar reasoning as in \cite[Lemma IV.2]{burger2014polyhedral}.
 Suppose by contradiction that there exist two processors $i$ and $p$ such that
  $\bar{J}^i>\bar{J}^p$. There exists $\eta_0>0$ such that
  $\bar{J}^i-\bar{J}^p>\eta_0$. Since the sequences $\{J(B^i(t))\}_{t>0}$ and
  $\{J(B^p(t))\}_{t>0}$ are monotonically increasing and convergent, for any
  $\eta>0$, there exists a time $T_\eta$ such that for all $t\geq T_\eta$,
  $\bar{J}^i-J(B^i(t))\leq\eta$ and $\bar{J}^p-J(B^p(t))\leq \eta$. This implies
  that there exists a $T_{\eta_0}$ such that for all $t\geq T_{\eta_0}$,
  $J(B^i(t))\geq\bar{J}^i-\eta_0>\bar{J}^p$. Additionally, since the objective
  function is increasing, it follows that for any time instant $t^\prime\geq 0$,
  $J(B^p(t^\prime))\leq \bar{J}^p$. Thus, for all $t\geq T_{\eta_0}$ and all
  $t^\prime\geq0$
\begin{equation}\label{eq: first contradiciton}
J(B^i(t))>J(B^p(t^\prime)).
\end{equation}
On the other hand, since the graph is uniformly jointly strongly connected,
there exists a time-varying path of length at most $nL$ from node $i$ to node
$p$. Therefore, for all $t\geq T_{\eta_0}$
\begin{equation}\label{eq: second contradiction}
J(B^i(t))\leq J(B^p(t+nL)).
\end{equation}
However, \eqref{eq: second contradiction} contradicts \eqref{eq: first
  contradiciton} proving that $J^1=\ldots=J^n$. Therefore, it must hold that
$\lim_{t\rightarrow\infty}|J(B^i(t))-J(B^j(t))|\rightarrow 0$ for all
$i,j\in\mathcal{V}$. This proves the second statement of the theorem.

%The third statement is proven by combining the second statement and  Assumption \ref{assum: nondegeneracy}.  
%Based on Assumption \ref{assum: nondegeneracy} any subproblem of \eqref{eq: uncertain LP} has a unique solution. We also proved in the second statement that all nodes converge to the same cost asymptotically. Therefore, all candidate solutions converge to a common point asymptotically and based on Assumption \ref{assum: nondegeneracy}, there are only $d$ constraints intersecting at the minimum point which form the basis. This proves the third statement of the theorem. 

% \item[2.] 
If $J(B^i(t_0))=J(B^i(t_0+2nL))$ and
considering the point that node $p$ can be any node of the graph, then all
nodes have the same cost. That is, $J(B^1(t))=\ldots=J(B^n(t))$. This combined with Assumption~\ref{assum: nondegeneracy} proves the third statement of the theorem. That is, if the candidate  solution is not updated for $2nL+1$ communication rounds, all nodes have a common  solution and hence the distributed algorithm can be halted.

% \end{enumerate}

%\MC{Proving that $\bar{J}=J^*$\\}
%In the previous part, we proved that all the nodes converge to the same objective value asymptoticly. Denote by $\bar{J}$ the objective value that all nodes converge to, that is $\bar{J}^1=\ldots=\bar{J}^n=\bar{J}$ also let $\bar{B}$ and $\bar{\theta}$ denote the basis and candidate solution corresponding to $\bar{J}$ \MC{I use the non-degeneracy assumption.}  Here, we prove that $\bar{J}=J^*$ almost surely. We proceed by contradiction. Let $\bar{J}<J^*$ \MC{with non-zero probability}. This implies that \MC{with non-zero probability} there exists a subset of uncertainty set $Q_v\subset\mathbb{Q}$ for which $J(\bar{B}\cup H(q_v))>\bar{J},q_v\in Q_v$ \MC{with non-zero probability}. Therefore, the candidate solution $\bar{\theta}$ is not robustly feasible \MC{with non-zero probability. }
%\begin{equation}\label{eq: proof termination point}
%\Pr\{q\in\mathbb{Q}:\bar{\theta}\notin\bigcap_{i=1}^n\mathcal{H}^i(q)\}>0.
%\end{equation}
%
%However, due to \MC{Assumption 000}, there is a non-zero probability of detecting a violation sample $q_\texttt{viol}$ at each step of Algorithm~\ref{alg: distributed algorithm}. Thus, with probability one Algorithm~\ref{alg: distributed algorithm} can not terminate at a point not being robustly feasible. In other words, the probability of terminating at a point not being robustly feasible is zero. This contradicts with \eqref{eq: proof termination point}. Therefore, $\bar{J}=J^*$ with probability one.
%

\textbf{Proof of forth statement:\\}
We first note that using  \cite[Theorem 1]{Chamanbaz_TAC_2016}, \cite[Theorem 3]{calafiore2007probabilistic} and \cite[Theorem 5.3]{dabbene2010randomized} we can show that---at any iteration $t$---if the sample size is selected based on \eqref{eq:sample bound Mk} and the verification step is successful, that is  $q_\texttt{viol}=\emptyset$, then 
\begin{equation*}
\Pr^M\left\{\ms\in\mathbb{Q}^M:\Pr\{q\in\mathbb{Q}:\theta^i(t)\notin\mathcal{H}^i(q)\}\leq\varepsilon_i \right\}\geq1-\delta_i.
\end{equation*}
We note that the above inequality is a centralized result and holds only for the agent's own constraint $\mathcal{H}^i(q)$. 
Also since for $\theta_\text{sol}$, the verification has  to be successful, then 
\begin{equation}\label{eq: proof misclass for one agent}
\Pr^M\left\{\ms\in\mathbb{Q}^M:\Pr\{q\in\mathbb{Q}:\theta_\text{sol}\notin\mathcal{H}^i(q)\}\leq\varepsilon_i \right\}\geq1-\delta_i.
\end{equation}
We further remark that the sample bound~\eqref{eq:sample bound Mk} is obtained by replacing $k_t-1$, $\delta/2$ and $\gamma$  in~\cite[Eq.~(10)]{Chamanbaz_TAC_2016} with $\infty$, $\delta_i$ and $1.1$ respectively, see~\cite[Remark 1]{Chamanbaz_TAC_2016} for a discussion on optimal value of $\gamma$. 

Now, we are interested in bounding the probability by which $\theta_{\text{sol}}\notin\bigcap_{i=1}^n\mathcal{H}^i(q)$, i.e. 
\begin{equation}\label{eq: proof misclass for multiple agent}
\Pr^M\left\{\ms\in\mathbb{Q}^M:\Pr\left\{q\in\mathbb{Q}:\theta_{\text{sol}}\notin\bigcap_{i=1}^n\mathcal{H}^i(q)\right\}\leq\varepsilon \right\}.
\end{equation}
In order to bound \eqref{eq: proof misclass for multiple agent}, we follow similar reasoning stated in \cite{margellos2016distributed}. 
Define the following events
\begin{align*}
&\Bad_i  \doteq\{\theta_{\text{sol}}\notin\mathcal{H}^i(q),\forall q\in\mathbb{Q}\}\\
&\Bad  \doteq\{\theta_{\text{sol}}\notin\bigcap_{i=1}^n\mathcal{H}^i(q),\forall q\in\mathbb{Q}\}.
\end{align*}
Equations \eqref{eq: proof misclass for one agent} and \eqref{eq: proof misclass for multiple agent} can be written in terms of the events \Bad$_i$ and \Bad
\begin{align}\label{eq: proof misclass based on Badi}
&\Pr^M\left\{\ms\in\mathbb{Q}^M:\Pr\{\Bad_i\}\leq\varepsilon_i\right\}\geq1-\delta_i\\ \label{eq: proof misclass based on Bad}
&\Pr^M\left\{\ms\in\mathbb{Q}^M:\Pr\{\Bad\}\leq\varepsilon\right\}
\end{align}
respectively. One can observe that 
\[
\theta_{\text{sol}}\notin\bigcap_{i=1}^n\mathcal{H}^i(q) \Rightarrow \exists i\in\{1,\ldots,n\}: \theta_{\text{sol}}\notin\mathcal{H}^i(q).
\]
Hence, the event $\Bad$  can be written as  the union of events $\Bad_i,~i=1,\ldots,n$, that is  $\Bad=\Bad_1\cup\Bad_2\cup\ldots\cup\Bad_n$.
Invoking Boole's inequality \citep{Comtet1974} (also known as Bonferroni's inequality), we have 
\begin{equation}\label{eq: Boole inequality}
\Pr\{\Bad\}\leq\sum_{i=1}^{n}\Pr\{\Bad_i\}.
\end{equation} 
Replacing $\Pr\{\Bad\}$ in \eqref{eq: proof misclass based on Bad} with the right hand side of \eqref{eq: Boole inequality} we obtain
\begin{align*}
 & \Pr^M\left\{\ms\in\mathbb{Q}^M:\sum_{i=1}^{n}\Pr\{\Bad_i\}\leq\varepsilon\right\}  \\
= &\Pr^M\left\{\ms\in\mathbb{Q}^M:\sum_{i=1}^{n}\Pr\{\Bad_i\}\leq\sum_{i=1}^{n}\varepsilon_i\right\}\\
\geq& \Pr^M\left\{\bigcap_{i=1}^n\left\{\ms\in\mathbb{Q}^M:\Pr\{\Bad_i\}\leq\varepsilon_i\right\}\right\}\\
\geq& 1- \sum_{i=1}^{n}\Pr^M\left\{\ms\in\mathbb{Q}^M:\Pr\{\Bad_i\}>\varepsilon_i\right\}\\
\geq& 1-\sum_{i=1}^{n}\delta_i= 1-\delta.
\end{align*}
We remark that the third line of above equation comes from the fact that if $\Pr\{\Bad_i\}\leq\varepsilon_i, ~\forall i=1,\ldots,n$ then, one can ensure that $\sum_{i=1}^{n}\Pr\{\Bad_i\}\leq\sum_{i=1}^{n}\varepsilon_i$. The forth line also is due to the fact that $\Pr\{\bigcap_iA_i\} = 1-\Pr\{\bigcup_i A_i^c\}$ where $A_i^c$ is the complement of the event $A_i$.

\textbf{Proof of fifth statement:\\}
We first note that if the solution $\theta_\text{sol}$ is violated for a sample $q_v$ from the set of uncertainty, that is, $\theta_\text{sol}\notin\bigcap_{i=1}^n\mathcal{H}^i(q_v)$, then $J(B_\text{sol}\cup H(q_v))\geq J(B_\text{sol})$ with $H(q_v)\doteq\bigcup_{i=1}^nH^i(q_v)$. However, due to Assumption~\ref{assum: nondegeneracy},  any subproblem of \eqref{eq: uncertain LP} has a unique minimum point and hence, $J(B_\text{sol}\cup H(q_v))\neq J(B_\text{sol})$. This argument combined with the result of forth statement proves the fifth statement of the theorem. That is, the probability that the solution $\theta_\text{sol}$  is no longer optimal for a new sample  equals the probability that the solution is violated by the new sample.  
\begin{flushright}

$\blacksquare$
\end{flushright}

%\begin{remark}[Probabilistic robustness and optimality]~\\
%Based on the forth statement of Theorem~\ref{thm: convergence }, the probability that the solution $\theta_\text{sol}$ is violated for a new sample (constraint) is smaller than $\varepsilon$ and this statement holds with probability at least $1-\delta$. Similarly, based on the fifth statement, the probability that the optimal objective $J(B_\text{sol})$ is no longer optimal for a new sample (constraint)  is smaller than $\varepsilon$ and this statement holds with probability at least $1-\delta$.
%\end{remark}

%-----------------------------------------------------------------------------------------

\section{Numerical Simulation}\label{sec: numerical simulations}
We test the effectiveness of the distributed algorithm presented in
Section~\ref{sec: distributed algorithm} through extensive numerical
simulations.  To this end, we generate random linear programs (LP)---with a
large number of uncertain parameters---assigned to various nodes of the
network. Each node is assigned an uncertain set of the form 
% has an uncertain LP of the form
\begin{align*}% \nonumber
% \min&~~c^T\theta\\ 
% \text{subject to}&~~ 
(A^0+A_q)^T\theta\leq b, %\label{eq: LP in Example}
\end{align*}
where $A^0$ is a fixed (nominal) matrix and $A_q$ is an interval matrix---a
matrix whose entries are bounded in given intervals---defining the uncertainty
in the optimization problem. We follow the methodology presented in
\cite{dunham_experimental_1977} in order to generate $A^0,~b$ and the problem
cost $c$ such that the linear program % \eqref{eq: LP in Example}
is always feasible. In particular, elements of $A^0$ are drawn from standard
Gaussian distribution (mean$=0$ and variance$=1$). The $i$-th element of $b$ is
define by $b_i = (\sum_{j=1}^{d}A^0_{ij})^{1/2}$. The objective direction
$c$---which is the same for all the nodes---is also drawn from the standard
Gaussian distribution. The communication graph $\mathcal{G}$ is a random
connected graph with fixed number of neighbors.
%In each node, denote by $A^0$ the nominal (not uncertain) matrix defining the constraint set of the agent. The   We consider a scenario where the matrix $A$ which defines the constraint set of the agent is uncertain. More precisely,  $A = A^0+A_q$  with $A_q$ being an interval matrix, i.e. a matrix whose entries are bounded in given intervals. Here, entries of $A_q$ are considered to be bounded within $[-\rho,\rho]$.
 \begin{table*}[t]
 \caption{The average---over all nodes---number of times a basis is transmitted to the neighbors,  average number of times verification is performed ($k_i$ at the convergence) and empirical violation of the computed solution ($\theta_\text{sol}$) over $10,000$ random samples for different number of nodes and neighbors in each node. The simulation is performed $100$ times for each row and average results are reported.}
 \begin{center}
 \scalebox{0.9}{
 %\centering{}%
 \begin{tabular}{c|c|c|c||c|c|c}
 \toprule
\# Nodes  & \# Neighbors & Graph&  \#Constraints &  \# Transmissions &  $k_i$ at convergence & Empirical \tabularnewline
 $n$  &  in each node & diameter & in each node &       (averaged)    & (averaged)    & violation \tabularnewline
 \midrule
 \midrule
 $10$  &  $3$  & $4$ &  $100$  & $29.57$  & $31.69$  & $2.81\times 10^{-4}$   \tabularnewline
 \midrule
 $20$  &  $4$  &  $4$ & $100$  & $26.92$   & $29.02$ & $1.7\times 10^{-4}$       \tabularnewline
 \midrule
$50$  &  $6$  &  $4$ & $100$  & $26.47$   & $28.51$ & $7\times 10^{-5}$       \tabularnewline
\midrule
$100$  &  $7$  &  $4$ & $100$  & $26.63$   & $28.68$ & $2.9\times 10^{-5}$       \tabularnewline
 \bottomrule
 \end{tabular}
 }
 \end{center}

 \label{tab: simulation results}
 \end{table*}
% \begin{figure}
% \begin{center}
% \includegraphics[width=8.4cm]{graph.eps}    % The printed column width is 8.4 cm.
% \caption{A random connected graph } 
% \label{fig:graph}
% \end{center}
% \end{figure}

 A workstation with $12$ cores and $48$ GB of RAM is used to emulate the network
 model. From an implementation viewpoint, each node executes Algorithm~\ref{alg:
   distributed algorithm} in an independent Matlab environment and the
 communication is modeled by sharing files between different Matlab
 environments.  We use the \texttt{linprog} function of Mosek
 \citep{andersen2000mosek} to solve optimization problems appearing at each
 iteration of the distributed algorithm.
%In Table~\ref{tab: simulation results} we vary different parameters such as uncertainty radius $\rho$, number of nodes $n$, number of neighbors of each node, number of constraints in each node $m$ and number of design variables $d$ and solve the resulting uncertain  distributed linear program.  

 In Table~\ref{tab: simulation results}, we change the number of nodes and
 neighbors such that the graph diameter is always $4$. The number of constraints
 in each node is also kept at $100$. We set the dimension of decision variables
 to $d=5$ and consider all elements of $A_q$ to be bounded in $[-0.2,0.2]$. The
 underlying probability distribution is selected to be uniform due to its
 worst-case nature.  The probabilistic accuracy and
 confidence levels of each agent ($\varepsilon_i$ and $\delta_i$) are $0.1/n$
 and $10^{-8}/n$ respectively, with $n$ being the number of nodes (first column
 of Table~\ref{tab: simulation results}). We report the average---over all
 nodes---number of times each node updates its basis and transmits it to the
 outgoing neighbors. It is assumed that each node keeps the latest information
 received from neighbors and hence, if the basis is not updated, there is no
 need to re-transmit it to the neighbors. This also accounts for the
 asynchronicity of the distributed algorithm.  We also report the average---over
 all nodes---number of times node $i$ performs the verification step, i.e.,
 $k_i$ at convergence. This allows us to show that with a small number of
 ``design'' samples used in the optimization step, nodes compute a solution with
 high degree of robustness.  In order to examine robustness of the obtained
 solution, we run an \emph{aposteriori} analysis based on Monte Carlo
 simulation. To this end, we collect all the constraints across different nodes
 in a single problem of form \eqref{eq: uncertain LP} and check---in a
 centralized setup---the feasibility of the obtained solution for $10,000$
 random samples extracted from the uncertain set. The empirical violation (last
 column of Table~\ref{tab: simulation results}) is measured by dividing the
 number of samples that violate the solution by $10,000$. Since
 Algorithm~\ref{alg: distributed algorithm} has a stochastic nature, we run the
 simulation $100$ times for each row of Table~\ref{tab: simulation results} and
 report the average values.

\begin{figure}
\begin{center}
\includegraphics[width=0.48\textwidth]{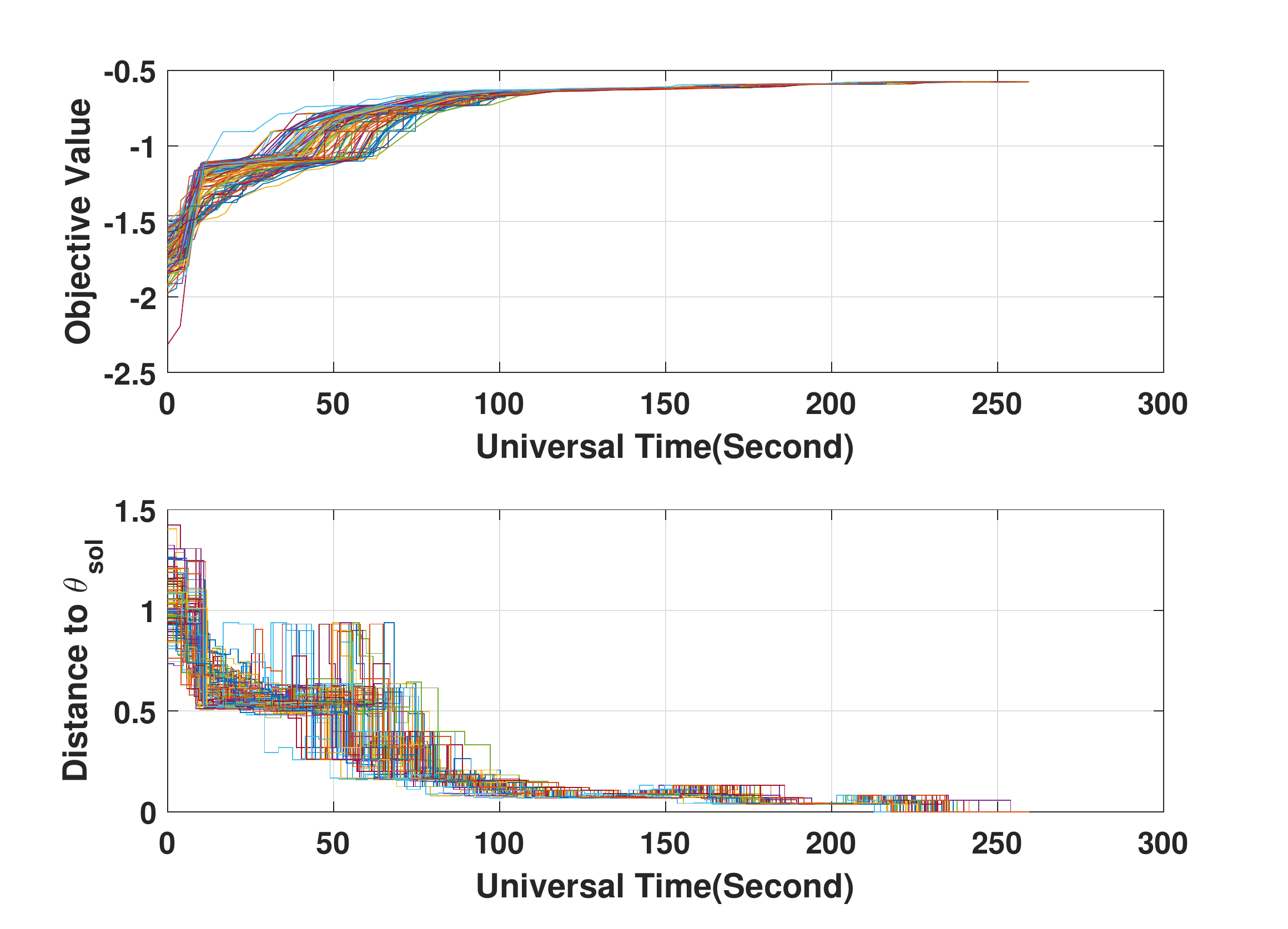}    % The printed column width is 8.4 cm.
\caption{Objective value and distance to $\theta_\text{sol}$ for all the nodes
  in the network corresponding to the last row of Table~\ref{tab: simulation
    results}.}  
\label{fig:objective and distance}
\end{center}
\end{figure}
 
We remark that it is very difficult in practice to check if
Assumption~\ref{assum: nondegeneracy} is satisfied for a given problem. However,
we observe that in all the simulations reported here, the candidate solutions
converge to the same point in finite time.  In Figure~\ref{fig:objective and
  distance}, we report the objective value and the distance of candidate
solutions $\theta^i(t), \forall i\in\{1,\ldots,n\}$ from $\theta_\text{sol}$
along the distributed algorithm execution for a problem instance corresponding
to the last row of Table~\ref{tab: simulation results}. It is observed that all
the nodes converge to the same solution $\theta_\text{sol}$.

%-----------------------------------------------------------------------------------------
\section{Conclusions}\label{sec: conclusions}
In this paper, we proposed a randomized distributed algorithm for solving robust
linear programs (LP) in which the constraint sets are scattered across a network
of processors communicating according to a directed time-varying graph. The
distributed algorithm has a sequential nature consisting of two main steps:
verification and optimization. Each processor iteratively verifies a candidate
solution through a Monte Carlo algorithm, and solves a local LP whose constraint
set includes its current basis, the collection of bases from neighbors and
possibly, a constraint---provided by the Monte Carlo algorithm---violating the
candidate solution. The two steps, i.e. verification and optimization, are
repeated till a local stopping criteria is met and all nodes converge to a
common solution. We analyze the convergence properties of the proposed
algorithm.

\bibliography{ifacconf}

\end{document}